\documentclass[12pt]{article}
\UseRawInputEncoding
\usepackage{amsfonts}
\usepackage{latexsym,amsmath}
\usepackage{amssymb,array}
\usepackage{setspace}
\makeatletter \oddsidemargin 0in \evensidemargin 0in \textwidth 16cm
 \RequirePackage[dvips]{graphicx} \textheight 18cm
\newcommand{\bt}{\begin{Theorem}}
\newcommand{\et}{\end{Theorem}}
\newcommand{\bi}{\begin{itemize}}
\newcommand{\ei}{\end{itemize}}
\newcommand{\bea}{\begin{eqnarray}}
\newcommand{\eea}{\end{eqnarray}}
\newtheorem{Definition}{Definition}[section]
\newtheorem{Theorem}[Definition]{Theorem}
\newtheorem{Lemma}[Definition]{Lemma}
\newcommand{\bl}{\begin{Lemma}}
\newcommand{\el}{\end{Lemma}}
\newcommand{\bpr}{\begin{Proposition}}
\newcommand{\epr}{\end{Proposition}}
\newtheorem{Notation}[Definition]{Notation}
\newcommand{\bn}{\begin{Notation}}
\newcommand{\en}{\end{Notation}}
\newtheorem{Remark}[Definition]{Remark}
\newcommand{\br}{\begin{Remark}}
\newcommand{\er}{\end{Remark}}
\newtheorem{Observation}[Definition]{Observation}
\newcommand{\bo}{\begin{Observation}}
\newcommand{\eo}{\end{Observation}}
\newtheorem{Construction}[Definition]{Construction}
\newcommand{\bcon}{\begin{Construction}}
\newcommand{\econ}{\end{Construction}}
\newtheorem{Example}[Definition]{Example}
\newcommand{\bex}{\begin{Example}}
\newcommand{\eex}{\end{Example}}
\newcommand{\bd}{\begin{Definition}}
\newcommand{\ed}{\end{Definition}}
\newcommand{\bs}{\begin{subsection}}
\newcommand{\es}{\end{subsection}}
\newtheorem{Proposition}[Definition]{Proposition}
\newtheorem{Corollary}[Definition]{Corollary}
\newcommand{\bc}{\begin{Corollary}}
\newcommand{\ec}{\end{Corollary}}

\newcommand{\be}{\begin{equation}}
\newcommand{\ee}{\end{equation}}
\date{}

\begin{document}
\title{A multi-objective multi-item solid transportation problem with vehicle cost, volume and weight capacity under fuzzy environment}
\author{Mouhya B. Kar$^1$, Pradip Kundu$^2$,
Samarjit Kar$^3$\footnote{Corresponding author, Address: NIT
Durgapur, W.B.-713209, India,  e-mail: $kar_{-}s_{-}k$@yahoo.com,
Telephone:+919434453186 },
 Tandra Pal$^4$ \\
\textit{\small{$^1$Department of Computer Science and Engineering,
Heritage Institute of Technology,}}\\\textit{\small{Kolkata 700107,
India}}\\\textit{\small{$^2$Department of Mathematics, Birla Institute of Technology Mesra, Ranchi-835215, India}}\\
\textit{\small{$^3$Department of Mathematics, National  Institute of
Technology,  Durgapur-713209, India}}
\\\textit{\small{$^4$Department of Computer Science, National Institute of Technology, Durgapur 713209, India}}} \maketitle \textbf{Abstract}
\vspace{0.4cm}\\
Generally, in transportation problem, full vehicles (e.g., light
commercial vehicles, medium duty and heavy duty trucks, etc.) are to
be booked, and transportation cost of a vehicle has to be paid
irrespective of the fulfilment of the capacity of the vehicle.
Besides the transportation cost, total time that includes travel
time of a vehicle, loading and unloading times of products is also
an important issue. Also, instead of a single item, different types
of items may need to be transported from some sources to
destinations through different types of conveyances. The optimal
transportation policy may be affected by many other issues like
volume and weight of per unit of product, unavailability of
sufficient number of certain types of vehicles, etc. In this paper,
we formulate a multi-objective multi-item solid transportation
problem by addressing all these issues. The problem is formulated
with the transportation cost and time parameters as fuzzy variables.
Using credibility theory of fuzzy variables, a chance-constraint
programming model is formulated, and is then transformed into the
corresponding deterministic form. Finally numerical example is
provided to illustrate the problem.
\section{Introduction} As a generalization of classical transportation
problem (TP), the solid transportation problem (STP) has been
extended by considering some extra constraints along with source
constraints and destination constraints, e.g. constraints due to
various types of goods, limited conveyance capacities, etc. The STP
was first presented by \cite{hale} considering the conveyance
capacity constraints. Recently the STP has been studied by many
researchers by describing it with many models and methods under
different uncertain environments (random, fuzzy, rough, etc.). In
decision making problems like transportation, the possible values of
the system parameters can not be always exactly determined. Dealing
with different types of uncertainties in many practical problems is
still an emerging problem. Fuzzy set theory is one of the most
convenient and accepted tool to deal with uncertainty. Some recent
works on both theoretical and application of fuzzy sets theory are
in the field of group decision making
\cite{pliu1,pliu2,pliu3,pliu4,pliu5}, image processing
\cite{chen,hass,jeon}, neural network \cite{cli,wang}, fault
detection \cite{fei,wu}, etc.

Different STP's with associated fuzzy parameters are considered by
many researchers. Bit et al. \cite{bit} introduced the fuzzy
programming model for multi-objective STP. Li et al. \cite{LiI}
developped an improved genetic algorithm to solve multi-objective
STP with fuzzy numbers. Jim´enez and Verdegay \cite{jime} solved a
fuzzy STPs using evolutionary algorithm based solution method. Yang
and Feng \cite{yang1} constructed different goal programming models
to solve a bicriteria STP with stochastic parameters. Yang and Liu
\cite{yang2} constructed chance-constrained programming models to
solve fixed charge STP with fuzzy parameters. Ojha et al.
 \cite{ojha} investigated an entropy based multi-objective STP with
transportation costs and route-wise travel times as general fuzzy
numbers. Giri et al. \cite{giri} presented fixed charge multi-item
STP with all associated parameters represented as triangular fuzzy
numbers. Dalman et al. \cite{dalm} proposed an interval fuzzy
programming approach to solve a multi-objective multi-item STP.
Sinha et al. \cite{sinh} solved a bi-objective STP with interval
type-2 fuzzy numbers. Das et al. \cite{das} solved an STP with the
associated parameters as trapezoidal type-2 fuzzy numbers. Kundu et
al. \cite{kund2,kund3,kund4} investigated various STP models with
type-1 and type-2 fuzzy parameters. In multi-item transportation
system, sometime more than one item are transported from some
sources to some destinations through some conveyances. In many
real-world situations, it is observed that several objectives are to
be considered and optimized at the same time. Then the corresponding
problem becomes a multi-objective problem.
\par
Most of the papers (as mentioned above) discussed STPs by
considering total available capacities (space) of conveyances and so
that considering transportation cost of each unit of product
transported. However, for transportation systems where full vehicles
(e.g., light commercial vehicles, medium duty and heavy duty trucks,
rail coaches, etc.) are to be considered for transportation of
products, different types of issues appear in formulation of the
problem. Like transportation cost of a vehicle which is irrespective
of the fact whether the capacity of the vehicle is filed up or not;
volume and weight capacities of the vehicles; limitation of number
of certain types of vehicles, etc. Also, previous works mainly
considered only travel time of vehicles, but besides the travel time
of a vehicle, loading and unloading times of products are also
important which depend upon both the vehicle types and product
characteristics. In this paper, we present a multi-objective
multi-item solid transportation model by addressing all these
issues. The presented problem is formulated with transportation time
and cost parameters as fuzzy variables. The problem is described in
detail in Section 3.\par The rest of this paper is organized as
follows. Section 2 discusses some basic concepts about fuzzy
variable. Section 3 describes the problem and formulates the model
mathematically. In Section 4, the methodology for the solution of
the problem is presented. Section 5 describes two techniques to
solve multi-objective optimization problems, namely, fuzzy
programming technique and global criteria method. Section 6
illustrates the problem numerically, and finally the paper is
concluded in Section 7.
\section{Preliminaries}
A fuzzy variable \cite{nahi} is defined as a function from the
possibility space $(\Theta,\emph{p},Pos)$ to the set of real numbers
$\Re$ to describe fuzzy phenomena, where possibility measure $(Pos)$
\cite{zade2,wang} of a fuzzy event $\{\tilde{\xi}\in B\}$, $B\subset
\Re$ is defined as $Pos\{\tilde{\xi}\in B\}=\sup_{x\in
B}~~\mu_{\tilde{\xi}}(x),where $ $\mu_{\tilde{\xi}}(x)$ is the
possibility distribution of $\tilde{\xi}$.\par For normalized fuzzy
variable ($\sup_{x\in \Re}~\mu_{\tilde{\xi}}(x)$=1), necessity
measure $(Nec)$ is defined as $ Nec\{\tilde{\xi}\in
B\}=1-Pos\{\tilde{\xi}\in B^c\}=1-\sup_{x\in
B^c}~~\mu_{\tilde{\xi}}(x)$ and the credibility measure \cite{liu1}
of $\{\tilde{\xi}\in B\}$ is defined as $Cr\{\tilde{\xi}\in
B\}=\frac{1}{2}(Pos\{\tilde{\xi}\in B\}+Nec\{\tilde{\xi}\in
B\}).$\par\textbf{Optimistic and pessimistic value \cite{liu2}:} Let
$\tilde{\xi}$ be a fuzzy variable and $\alpha\in[0,1]$. Then,
$\alpha$-optimistic value and $\alpha$-pessimistic values of
$\tilde{\xi}$ are defined as follows.
$$\tilde{\xi}_{sup}(\alpha)=\sup\{r:Cr\{\tilde{\xi}\geq r\}\geq \alpha\},$$
$$\tilde{\xi}_{inf}(\alpha)=\inf\{r:Cr\{\tilde{\xi}\leq r\}\geq \alpha\}.$$ \begin{Example} Let
$\tilde{\xi}=(r_{1},r_{2},r_{3},r_{4})$ be a trapezoidal fuzzy
variable. Then its $\alpha$-optimistic and $\alpha$-pessimistic
values are as given below.
$$\tilde{\xi}_{sup}(\alpha)=\left\{
                                                                  \begin{array}{ll}
                                                                    2\alpha r_{3}+(1-2\alpha)r_{4}, & \hbox{if $\alpha\leq 0.5$;} \\
                                                                 (2\alpha-1)r_{1}+2(1-\alpha)r_{2}, & \hbox{if $\alpha>0.5$.}
                                                                  \end{array}
                                                                \right.$$  $$\tilde{\xi}_{inf}(\alpha)=\left\{
                                                                  \begin{array}{ll}
                                                                    (1-2\alpha) r_{1}+2\alpha r_{2}, & \hbox{if $\alpha\leq 0.5$;} \\
                                                                    2(1-\alpha) r_{3}+(2\alpha-1)r_{4}, & \hbox{if $\alpha>0.5$.}
                                                                  \end{array}
                                                                \right.$$
\end{Example}
\section{Problem description and model formulation}
In multi-item solid transportation problem (MISTP), different types
of items/products are to be transported from some sources to some
destinations through some conveyances (modes of transportation) so
that the objective (e.g., total transportation cost, time, profit,
etc.) is optimum. In many real transportation systems, full vehicles
(e.g., light commercial vehicles, medium duty and heavy duty trucks
for road transportation, coaches for rail transportation, etc.) are
to be booked and number of vehicles are determined according to the
amount of products to be transported through a particular route. In
such cases, full transportation cost of a vehicle has to be paid
irrespective of the capacity of the vehicle is filed up by the
products or not. So the allocation of the products are to be done in
such a way so that the volume capacities of the vehicles are field
up as much as possible. Time (transportation duration) is also an
important issue in transportation system. Beside travel time of a
vehicle; loading and unloading times of products are also important.
Loading and unloading times depend upon both the vehicle types and
product characteristics. In this paper, we have considered travel
time, and loading and unloading times of products for each type of
vehicles. We also consider the weight capacities of vehicles. Also
the number of vehicles of certain type of conveyance may be limited
for some route. In this situation, a constraint on the number of
available vehicles should be considered. This limitation of number
of vehicles can affect the optimal transportation policy. For
example, unavailability of sufficient number of vehicles of certain
type of conveyance may force to use another type of conveyance for
which cost is more.
\subsection{Model formulation}
Different parameters and decision variables as used to
formulate the mathematical model are given bellow:\\
\begin{tabbing}
aaaaaaaaaaaaaa \= abababababababababababab \kill
  Parameters \>  \\
  $p$ \>  Items/products; $p=1,2,...,l$\\
  $i$ \>  Source/origin; $i=1,2,...,m$\\
  $j$\> Destination/demand point; $j=1,2,...,n$\\
  $k$ \> Types of vehicles (modes of transportation); $k=1,2,...,K$\\
  $t_{ijk}$\> Time required to travel from source $i$ to destination $j$ through\\\> vehicle of type $k$\\
  $\alpha_{pk}$ \> Time of loading and unloading one unit of item $p$ for the vehicle of type k\\
  $c_{ijk}$ \> Per trip transportation cost of a vehicle of type $k$ for traveling from origin $i$\\ \> to destination $j$\\
  $V_{k}$ \> Volume capacity of a vehicle of type $k$\\
  $W_{k}$ \> Weight capacity of a vehicle of type $k$\\
  $v^p$ \> Volume of one unit of product $p$\\
  $w^p$ \> Weight of one unit of product $p$\\
  $a_{i}^p$ \> Amount of a product $p$ available at origin $i$\\
  $b_{j}^p$ \> Demand of the product $p$ at destination $j$\\
  $Q_{k}$ \> Number of available vehicles of type $k$\\
  $f_{i}$ \> The objective value
 \end{tabbing}
\begin{tabbing}
aaaaaaaaaaaaaa \= abababababababababababab \kill
  Decision variables \>  \\
  $x_{ijk}^p$ \>  Amount of item $p$ transported from source $i$
to destination $j$ using \\ \> vehicles of type $k$ \\
  $z_{ijk}$ \> Number of required vehicles of type $k$ for transporting goods from \\\> source $i$ to
destination $j$\\
   \end{tabbing}
The proposed bi-objective MISTP model with vehicle cost, volume and
weight capacity is formulated as follows:
\begin{eqnarray}\label{eqmistp1}\text{Min}~f_{1}&=&\sum_{i=1}^m\sum_{j=1}^n\sum_{k=1}^K~~c_{ijk}~z_{ijk},\\
\label{eqmistp2}Min~f_{2}&=&\sum_{i=1}^m\sum_{j=1}^n\sum_{k=1}^K~(t_{ijk}~z_{ijk}+\sum_{p=1}^l
\alpha_{pk}~x_{ijk}^p),\end{eqnarray}
\begin{eqnarray}\hspace{-2cm}\label{eqmistp3} \text{subject~to}~
\sum_{j=1}^n\sum_{k=1}^K x_{ijk}^p~&\leq& a_{i}^p,~~i=1,2,...,m;
p=1,2,...l,\\\label{eqmistp4} \sum_{i=1}^m\sum_{k=1}^K
x_{ijk}^p~~&\geq& b_{j}^p,~~j=1,2,...,n; p=1,2,...,l,
\\\label{eqmistp5} \sum_{p=1}^l v^p x_{ijk}^p&\leq& z_{ijk}\cdot
V_{k},~i=1,2,...,m; j=1,2,...,n; k=1,2,...,K,\\\label{eqmistp6}
\sum_{p=1}^l w^p x_{ijk}^p&\leq& z_{ijk}\cdot W_{k},~i=1,2,...,m;
j=1,2,...,n; k=1,2,...,K,\\\label{eqmistp7}
\sum_{i=1}^m\sum_{j=1}^n~z_{ijk}&\leq& Q_{k},~k=1,2,...,K,
\end{eqnarray}
\begin{equation}\label{eqmistp8}x_{ijk}^p\geq 0,~ \forall~ i,j,k,p.\end{equation}
Here, the first objective function, i.e., $f_{1}$ represents the
total transportation cost, and the second objective function, i.e.,
$f_{2}$ represents the total time (trip durations), where $t_{ijk}$
represents the the travel time for each vehicle of type $k$ from
source $i$ to destination $j$, and $\sum_{p=1}^l
\alpha_{pk}~x_{ijk}^p$ represents the loading and unloading time of
all types of items transported from source $i$ to destination $j$
for the vehicle of type $k$.\par The constraint (\ref{eqmistp3})
ensures that total transported amount of each type of item from some
source must be equal to or less than the availability $(a_{i}^p)$ of
the item in that source. The constraint (\ref{eqmistp4}) indicates
that total transported amount of each type of item from the sources
should satisfy the demand $(b_{j}^p)$ of destination. The constraint
(\ref{eqmistp5}) ensures that total transported amount of products
must be equal to or less than the total volume capacities of all
types of allocated vehicles from a source $i$ to a destination $j$.
The constraint (\ref{eqmistp6}) ensures that weights of total
transported products must be equal to or less than the total weight
capacities of all types of allocated vehicles from a source $i$ to a
destination $j$. The constraint (\ref{eqmistp7}) is imposed on the
availability of vehicles of type $k$ for the source $i$ to
destination $j$.
\subsection{The MISTP Model with fuzzy transportation cost and time parameters}
Transportation cost depends upon fuel price, labor charge, tax
charge, etc., each of which fluctuate with time. So it is not easy
to predict the exact transportation cost of a route. Similarly,
travel time of vehicles depends upon condition of road, road jam,
vehicle condition; and loading and unloading times depend upon
availability of manpower, product characteristics, vehicle type,
etc. Generally, possible values of parameters are given by the
experts by approximate numbers, intervals, linguistic terms, etc.
Also each of the point in the given interval may not have the same
importance or possibility. For a large data set of a certain
parameter collected from previous experiments, generally all the
data are not equally possible. Such types of linguistic information,
approximate intervals, non equipossible data set can be expressed by
fuzzy numbers(/variables), especially by triangular or trapezoidal
fuzzy numbers \cite{bect,buck,dubo,kund2}.\par Consider that
transportation cost $c_{ijk}$, travel time $t_{ijk}$, loading and
unloading time $\alpha_{pk}$ in the above model are
represented by fuzzy variable respectively as follows:\\
$\tilde{c}_{ijk}=(c_{ijk}^{1},c_{ijk}^{2},c_{ijk}^{3},c_{ijk}^{4})$,
$\tilde{t}_{ijk}=(t_{ijk}^{1},t_{ijk}^{2},t_{ijk}^{3},t_{ijk}^{4})$,
$\tilde{\alpha}_{pk}=(\alpha_{pk}^1,\alpha_{pk}^2,\alpha_{pk}^3,\alpha_{pk}^4)$
for all $i,j,k$ and $p$. Then the problem (1)-(8) becomes
\begin{eqnarray}\label{fmistp1}&&\text{Min}~\tilde{f}_{1}=\sum_{i=1}^m\sum_{j=1}^n\sum_{k=1}^K~~\tilde{c}_{ijk}~z_{ijk},\\
\label{fmistp2}&&\text{Min}~\tilde{f}_{2}=\sum_{i=1}^m\sum_{j=1}^n\sum_{k=1}^K~(\tilde{t}_{ijk}~z_{ijk}+\sum_{p=1}^l
\tilde{\alpha}_{pk}~x_{ijk}^p),\\&& \label{conset}
\text{subject~to}~(\ref{eqmistp3})-(\ref{eqmistp8}).\end{eqnarray}
Since $\tilde{c}_{ijk}$ are trapezoidal fuzzy numbers and
$z_{ijk}\geq 0$ for all $i,j,k,$ so
$\tilde{f}_{1}=\sum_{i=1}^m\sum_{j=1}^n\sum_{k=1}^K~~\tilde{c}_{ijk}~z_{ijk}$
is also trapezoidal fuzzy number for any feasible solution and given
by $\tilde{f}_{1}=(r_{1},r_{2},r_{3},r_{4})$, where
\begin{equation}\label{eqr1r2} r_{1}=\sum_{i=1}^m\sum_{j=1}^n\sum_{k=1}^K~~c_{ijk}^{1}~z_{ijk}~
,~r_{2}=\sum_{i=1}^m\sum_{j=1}^n\sum_{k=1}^K~~c_{ijk}^{2}~z_{ijk},\end{equation}
\begin{equation}\label{eqr3r4} r_{3}=\sum_{i=1}^m\sum_{j=1}^n\sum_{k=1}^K~~c_{ijk}^{3}~z_{ijk}~
,~r_{4}=\sum_{i=1}^m\sum_{j=1}^n\sum_{k=1}^K~~c_{ijk}^{4}~z_{ijk}.\end{equation}
Similarly $\tilde{t}_{ijk}$, $\tilde{\alpha}_{pk}$ are trapezoidal
fuzzy numbers and $x_{ijk}^p\geq 0$ for all $i,j,k,p$. So
$\tilde{f}_{2}$ can be represented by
$\tilde{f}_{2}=(s_{1},s_{2},s_{3},s_{4})$, where
\begin{equation}\label{eqs1s2} s_{1}=\sum_{i=1}^m\sum_{j=1}^n\sum_{k=1}^K~(t_{ijk}^1~z_{ijk}+\sum_{p=1}^l
\alpha_{pk}^1~x_{ijk}^p)~
,~s_{2}=\sum_{i=1}^m\sum_{j=1}^n\sum_{k=1}^K~(t_{ijk}^2~z_{ijk}+\sum_{p=1}^l
\alpha_{pk}^2~x_{ijk}^p),\end{equation}
\begin{equation}\label{eqs3s4} s_{3}=\sum_{i=1}^m\sum_{j=1}^n\sum_{k=1}^K~(t_{ijk}^3~z_{ijk}+\sum_{p=1}^l
\alpha_{pk}^3~x_{ijk}^p)~
,~s_{4}=\sum_{i=1}^m\sum_{j=1}^n\sum_{k=1}^K~(t_{ijk}^4~z_{ijk}+\sum_{p=1}^l
\alpha_{pk}^4~x_{ijk}^p).\end{equation} \section{Solution
methodology: Chance-constrained programming} Chance-constrained
programming with fuzzy parameters was developed by Liu and Iwamura
\cite{liui}, Liu \cite{liu}, Yang and Liu \cite{yang2}, Kundu et al.
\cite{kund2} and many more authors. This method is used to solve the
problems with chance-constraints. In this method, the uncertain
constraints are allowed to be violated such that constraints must be
satisfied at some chance (confidence) level. Applying this method
using credibility measure for the above problem (given in Section
3.2) with fuzzy transportation costs and time parameters, the
following chance-constrained programming (CCP) model is formulated:
\begin{equation}\label{minf1bar}\text{Min~(Min}_{\bar{f_{1}}}~~\bar{f_{1}}),\end{equation}
\begin{equation}\label{minf2bar}\text{Min~(Min}_{\bar{f_{2}}}~~\bar{f_{2}}),\end{equation}
\begin{equation}\label{crf1} s.t.~Cr\left\{\sum_{i=1}^m\sum_{j=1}^n\sum_{k=1}^K~~\tilde{c}_{ijk}~z_{ijk}\leq
\bar{f_{1}}\right\}\geq \eta,\end{equation}
\begin{equation}\label{crf2} Cr\left\{\sum_{i=1}^m\sum_{j=1}^n\sum_{k=1}^K~(\tilde{t}_{ijk}~y_{ijk}+\sum_{p=1}^l
\tilde{\alpha}_{pk}~x_{ijk}^p)\leq \bar{f_{2}}\right\}\geq
\gamma,\end{equation} \begin{equation} \text{subject~to}~
(\ref{eqmistp3})-(\ref{eqmistp8}).\end{equation} Since our problem
is minimization problem, for the objective functions (\ref{fmistp1})
and (\ref{fmistp2}) we want to minimize $\eta$-pessimistic and
$\gamma$-pessimistic values of $\tilde{f}_1$ and $\tilde{f}_2$
respectively, where $\eta$ and $\gamma$ $(0<\eta,\gamma\leq 1)$ are
preassigned values. More specifically, for the objective function
(\ref{fmistp1}) we want to minimize
$\inf\{\bar{f_{1}}:Cr\{\tilde{f}_1\leq \bar{f_{1}}\}\geq \eta\}$
which is represented by (\ref{minf1bar}) and (\ref{crf1}) together.
Similar explanation follows for (\ref{minf2bar}) and (\ref{crf2}).
%so that (\ref{minf1bar}) and (\ref{crf1}) together represent that we
%are going to minimize smallest possible $\bar{f_{1}}$ with the
%condition that the credibility degree that the respective objective
%function $\tilde{f}_1$ less than or equal to it should be at least
%the preassigned numbers $\eta$.
\subsection{Deterministic form of the CCP Model}
In the above CCP model, $\text{Min}_{\bar{f_{1}}}~\bar{f_{1}}$, s.t.
$Cr\{\sum_{i=1}^m\sum_{j=1}^n\sum_{k=1}^K~\tilde{c}_{ijk}~z_{ijk}\leq
\bar{f_{1}}\}\geq \eta$ can be equivalently computed as
$f'_{1}=\inf\{r:Cr\{\tilde{f}_{1}\leq r\}\geq \eta\}$, which is
nothing but $\eta$-pessimistic value to $\tilde{f}_{1}$
$(i.e.,~\tilde{f}_{1_{inf}}(\eta))$ and so is equal to $f'_{1}$
,where
\begin{equation}\label{crispf1} f'_{1}=\left\{
\begin{array}{ll}
(1-2\eta)r_{1}+2\eta~r_{2}, & \hbox{if~$\eta\leq 0.5$;} \\
2(1-\eta)r_{3}+(2\eta-1)r_{4}, & \hbox{if~$\eta>0.5$.}
\end{array}
\right.\end{equation} Here $r_{1}$, $r_{2}$, $r_{3}$, $r_{4}$ are
given in equations (\ref{eqr1r2}) and (\ref{eqr3r4}).\par Similarly
$\text{Min}_{\bar{f_{2}}}~~\bar{f_{2}},$ s.t.
$Cr\{\sum_{i=1}^m\sum_{j=1}^n\sum_{k=1}^K~(\tilde{t}_{ijk}~y_{ijk}+\sum_{p=1}^l
\tilde{\alpha}_{pk}~x_{ijk}^p)\leq \bar{f_{2}}\}\geq \gamma$ is
equivalent to $f'_{2}=\inf\{s:Cr\{\tilde{f}_{2}\leq s\}\geq
\gamma\}$, which is nothing but $\gamma$-pessimistic value to
$\tilde{f}_{2}$ $(i.e.~\tilde{f}_{2_{inf}}(\gamma))$ and so is equal
to $f'_{2}$ , where
\begin{equation}\label{crispf2} f'_{2}=\left\{
\begin{array}{ll}
(1-2\gamma)s_{1}+2\gamma~s_{2}, & \hbox{if~$\gamma\leq 0.5$;} \\
2(1-\gamma)s_{3}+(2\gamma-1)s_{4}, & \hbox{if~$\gamma>0.5$.}
\end{array}
\right.\end{equation} Here $s_{1}$, $s_{2}$, $s_{3}$, $s_{4}$ are
given in equations (\ref{eqs1s2}) and (\ref{eqs3s4}).\par Finally
crisp form of the above CCP Model can be written as
\begin{eqnarray*}&&\text{Min}~f'_{1},\\&& \text{Min}~f'_{2},\\&&\text{subject~to}~
(\ref{eqmistp3})-(\ref{eqmistp8}).
\end{eqnarray*}
\par To solve the deterministic multi-objective problem, we apply two
multi-objective optimization methods, namely, the fuzzy programming
technique \cite{zimm,bit} and global criterion method which are
discussed briefly in the next section.
\section{Techniques used to solve multi-objective optimization problem}
Consider a multi-objective optimization problem with $R$ objective
functions: \begin{eqnarray*}\text{Min}~F(x)&=&\left(f_1(x),f_2(x),...,f_R(x)\right)^T\\
s.t.~x&\in& D,
\end{eqnarray*}
where $D$ is the set of feasible solutions.
\subsection{Fuzzy Programming Technique} Zimmermann \cite{zimm} first introduced fuzzy linear programming approach for solving problem with multiple
objectives and he showed that fuzzy linear programming always gives
efficient solutions and an optimal compromise solution. The steps to
solve the multi-objective models using fuzzy programming technique
are as follows:
\\Step 1: Solve the multi-objective problem as a single objective
problem using, each time, only one objective $f_{t}~(t=1,2,...,R)$
(ignore all other objectives) to obtain the optimal solution
$x^{t*}$ of R different single objective problem.\\ Step 2:
Calculate the values of each objective function at all these R
optimal solutions $x^{t*}~(t=1,2,...,R)$ and find the upper and
lower bound for each objective given by
$U_{t}=Max\{f_{t}(x^{1*}),f_{t}(x^{2*}),...,f_{t}(x^{R*})\}$ and $L_{t}=f_{t}(x^{t*})$, $t=1,2,...,R$ respectively.\\
Step 3: Then an initial model is given by $$Find~~x$$
$$\text{subject~to} ~~~~f_{t}(x)\leqslant L_{t},~~~t=1,2,...,R$$
\hspace{6.5cm} and $x \in D$.\\ However, generally due to
conflicting nature of the objective functions, feasible solution of
the above model does not always exists.\\ Step 4: Construct the
linear membership function $\mu_{t}(f_{t})$ corresponding to $t$-th
objective as
$$\mu_{t}(f_{t})=\left\{
                                      \begin{array}{ll}
                                        1, & \hbox{if $f_{t}\leq L_{t}$;} \\
                                        \frac{U_{t}-f_{t}(x)}{U_{t}-L_{t}}, & \hbox{if $L_{t}<f_{t}<U_{t}$;} \\
                                        0, & \hbox{if $f_{t}\geq U_{t}$,~~~~~~~~~~~~~~~~~~~~~~$\forall~ t$.}
                                      \end{array}
                                    \right.$$
Step 5: Formulate fuzzy linear programming problem using max-min
operator as  $$Max ~\lambda$$ $$subject ~to~\lambda\leq
\mu_{t}(f_{t})=(U_{t}-f_{t})/(U_{t}-L_{t}),~\forall t$$
\hspace{6.5cm} and $x \in D,$
$$\lambda\geq 0 ~and~ \lambda=min_{t}\{\mu_{t}(f_{t})\}.$$ Step-6: Now the reduced
problem is solved by a linear optimization technique and the optimum
compromise solutions are obtained.
\subsection{Global Criteria Method} Global criteria method gives a compromise solution for a multi-objective optimization problem. Actually this method
is a way of achieving compromise in minimizing the sum in deviations
of the ideal solutions (minimum value of the each objectives in case
of minimization problem) from the respective objective functions.
The steps of this method to solve the multi-objective models are as
follows:
\\Step-1: Solve the multi-objective problem as a single objective
problem using, each time, only one objective
$f_{t}~(t=1,2,...,R)$ ignoring all other objectives.\\
Step-2: From the results of step-1, determine the ideal objective
vector, say $(f_{1}^{\min},f_{2}^{\min},...,f_{R}^{\min})$ and
corresponding values of $(f_{1}^{\max},f_{2}^{\max},...,f_{R}^{\max})$.\\
Step-3: Formulate the following auxiliary problem $$Min~G(x)$$
\hspace{6.5cm} s.t. $x \in D,$
$$G(x)=Min\bigg\{\sum_{t=1}^R\bigg(\frac{f_{t}(x)-f_{t}^{\min}}{f_{t}^{\min}}\bigg)^q\bigg\}^{\frac{1}{q}},$$
$$or,~G(x)=Min\bigg\{\sum_{t=1}^R\bigg(\frac{f_{t}(x)-f_{t}^{\min}}{f_{t}^{\max}-f_{t}^{\min}}\bigg)^q\bigg\}^{\frac{1}{q}},$$
where $1\leq q\leq \infty$. An usual value of $q$ is 2. This method
is then called global criterion method in $L_{2}$ norm.
\section{Numerical Experiment}
To illustrate the MISTP model (\ref{fmistp1})-(\ref{conset}), we
consider a transportation plan in which two steel products
manufacturing company supply two types of steel products to three
cities by means of two types of conveyances which are super heavy
duty truck (dump truck) and heavy duty truck. That is, here $i=1,2$;
$j=1,2,3$; $k=1,2$ and $p=1,2$. Now, the problem is to make a
transportation plan for the next quarter such that the total
transportation cost and total transportation time are minimized at
the same time. To cope with uncertainty about the transportation
cost and transportation time, these parameters are considered as
trapezoidal fuzzy variables. The values of the remaining parameters
such as volume and weight capacity of each type of vehicle,
availability of product, demand of product, and maximum availability
of vehicles are deterministic. The transportation costs for two
types of vehicles for this problem are given in Table 1 and Table 2.
Travel time of vehicles are given in Table 3 and Table
4.\begin{table}[h]
\begin{center} \caption{Vehicle costs
$\tilde{c_{ij1}}$}
\begin{tabular}{|c|ccc|ccc|}
\hline $i\setminus j$ & \multicolumn{1}{c}{1} &
\multicolumn{1}{c}{2} & \multicolumn{1}{c|}{3} \\\hline $1$ &
(101,102,104,105) & (103,104,105,106) & (104,106,108,110)
 \\
$2$  & (102,104,106,107) & (108,110,111,112) & (102,103,104,106)
\\
\hline
\end{tabular}
\end{center}
\end{table}
\begin{table}[h]
\begin{center} \caption{Vehicle costs
$\tilde{c_{ij2}}$}
\begin{tabular}{|c|ccc|ccc|}
\hline $i\setminus j$ & \multicolumn{1}{c}{1} &
\multicolumn{1}{c}{2} & \multicolumn{1}{c|}{3} \\\hline $1$  &
(90,91,92,93) & (87,88,89,91) & (94,95,96,97)
 \\
$2$  & (94,96,97,98) & (92,93,94,96) & (93,94,95,97)\\ \hline
\end{tabular}
\end{center}
\end{table}
\begin{table}[h]
\begin{center} \caption{Travel time
$\tilde{t_{ij1}}$ (in hour)}
\begin{tabular}{|c|ccc|ccc|}
\hline $i\setminus j$ & \multicolumn{1}{c}{1} &
\multicolumn{1}{c}{2} & \multicolumn{1}{c|}{3} \\\hline $1$ &
(5,5.5,6,6.2) & (5.4,5.8,6,6.4) & (5.5,5.8,6,6.5)
 \\
$2$  & (5,5.5,6.2,6.4) & (5.8,6,6.5,6.8) & (5,5.5,6,6.4)
\\
\hline
\end{tabular}
\end{center}
\end{table}
\begin{table}[h]
\begin{center} \caption{Travel time
$\tilde{t_{ij2}}$ (in hour)}
\begin{tabular}{|c|ccc|ccc|}
\hline $i\setminus j$ & \multicolumn{1}{c}{1} &
\multicolumn{1}{c}{2} & \multicolumn{1}{c|}{3} \\\hline $1$ &
(4.6,5,5.5,5.6) & (4.5,4,8,5.4,5.6) & (4.8,5,5.2,5.8)
 \\
$2$ & (4.5,5,5.4,5.8) & (5,5.4,5.6,6) & (4.8,5,5.4,5.8)
\\
\hline
\end{tabular}
\end{center}
\end{table}The time of loading and unloading (in minute) of one unit of item $p$
into conveyance of type $k$ is $\tilde{\alpha}_{11}=(8,8.5,9,10)$,
$\tilde{\alpha}_{12}=(7.5,8,8.5,9)$,
$\tilde{\alpha}_{21}=(7,8,8.5,9)$ and
$\tilde{\alpha}_{22}=(6,7,8,8.5)$. Values of the remaining
parameters are given in Table 5.\par Now to solve the
\begin{table}[h]
\begin{center} \caption{Parameter values}
\begin{tabular}{|c|c|}
\hline Volume capacity of vehicles(in $ft^3$) & $V_{1}=406.12$,
$V_{2}=348$
\\\hline Weight capacity of vehicles(in kg) & $W_{1}=18400$, $W_{2}=15767$
 \\\hline
Volume of one unit of product(in $ft^3$) & $v^1=19.94$, $v^2=12.66$
\\
\hline Weight of one unit of product(in kg) &  $w^1=45$, $w^2=40$ \\
\hline Availability of product &  $a_{1}^1=625$, $a_{2}^1=428$, $a_{1}^2=450$, $a_{2}^2=380$ \\
\hline Demand of product &  $b_{1}^1=340$, $b_{2}^1=360$, $b_{3}^1=345$, $b_{1}^2=275$, $b_{2}^2=250$, $b_{3}^2=280$ \\
\hline Availability of vehicles & $Q_{1}=52$, $Q_{2}=35$\\ \hline
\end{tabular}
\end{center}
\end{table}
problem, we model it using CCP technique and consider the
credibility degrees $\eta=\gamma=0.9$. Then using equations
(\ref{crispf1}) and (\ref{crispf2}) we have the crisp form of the
proposed STP model as
\allowdisplaybreaks{\begin{eqnarray}\label{eqexstart}\text{Min}~f'_{1}&=&0.2~r_{3}+0.8~r_{4},\\
\text{Min}~f'_{2}&=&0.2~s_{3}+0.8~s_{4},\\\label{eqsub1}\hspace{-2cm}
\text{subject~to}~ \sum_{j=1}^3\sum_{k=1}^2 x_{ijk}^p~&\leq&
a_{i}^p,~~i=1,2; p=1,2,\\
\sum_{i=1}^2\sum_{k=1}^2 x_{ijk}^p~~&\geq&
b_{j}^p,~~j=1,2,3; p=1,2 \\
\sum_{p=1}^2 v^p x_{ijk}^p&\leq& z_{ijk}\cdot V_{k},~i=1,2,;
j=1,2,3; k=1,2,\\ \sum_{p=1}^2 w^p x_{ijk}^p&\leq&
z_{ijk}\cdot W_{k},~i=1,2,; j=1,2,3; k=1,2,\\
\sum_{i=1}^2\sum_{j=1}^3~z_{ijk}&\leq&
Q_{k},~k=1,2,\\\label{eqexend}x_{ijk}^p&\geq& 0,~ \forall~ i,j,k,p,
\end{eqnarray}}
where $r_{3}$, $r_{4}$ and $s_{3}$, $s_{4}$ are given in
\ref{eqr3r4} and \ref{eqs3s4}.
\par To solve the multi-objective problem described in (\ref{eqexstart})-(\ref{eqexend}), we use two multi-objective optimization
methods.\par Based on the fuzzy programming technique (c.f. Sec.
5.1) with auxiliary variable $\lambda$, the STP problem
(\ref{eqexstart})-(\ref{eqexend}) can be modeled as follows in
(\ref{eqflp})-(\ref{eqflpend}).
\begin{eqnarray}\label{eqflp}\text{Min}~\lambda,\\\text{subject~to}~
f'_{1}+\lambda(8211.6-8166.6)&\leq&8211.6,\\
f'_{2}+\lambda(785.95-770.1767)&\leq&785.95,\\\label{eqflpend}\text{and}~(\ref{eqsub1})-(\ref{eqexend}).
\end{eqnarray}
\begin{table}[h]\label{tabresultfp}
\begin{center} {Table 6: Compromise solution using fuzzy programming
technique
\begin{tabular}{|l|}
               \hline
               $z_{111} = 13$, $z_{121} = 5$, $z_{211} = 8$, $z_{231} = 24$, $z_{112} = 5$, $z_{122}
               =24$, $z_{132}= 1$, $z_{232} = 1$,\\ $x_{111}^1 = 153$, $x_{111}^2 =
176$, $x_{121}^1 = 1$, $x_{121}^2 =156$, $x_{211}^1 = 100$,
$x_{211}^2 = 99$, $x_{231}^1 = 312$,\\$x_{231}^2 = 278$, $x_{112}^1
= 87$, $x_{122}^1 = 359$, $x_{122}^2 = 94$,
$x_{132}^1 = 17$, $x_{232}^1 = 16$, $x_{232}^2 = 2$,\\$\lambda=0.7077$, $Min~f'_{1}=8177.4$, $Min~f'_{2}=774.7867$. \\
\hline
\end{tabular} }
\end{center}
\end{table}
Solving the problem (\ref{eqflp})-(\ref{eqflpend}), we have the
compromise solution for the multi-objective STP defined in
(\ref{eqexstart})-(\ref{eqexend}) as presented in Table 6.\par
Solution using global criterion method: \par Based on the global
criterion method in $L_{2}$ norm (c.f. Sec. 5.2), we have the
following problem
\begin{eqnarray}\label{eqglc}&&\text{Min}~G=\left\{\left(\frac{f'_{1}-8166.6}{8166.6}\right)^2 +\left(\frac{f'_{2}-770.1767}{770.1767}\right)^2\right\}^{\frac{1}{2}},\\&&\text{subject~to}~
\label{eqglce}(\ref{eqsub1})-(\ref{eqexend}).
\end{eqnarray}
Solving the problem (\ref{eqglc})-(\ref{eqglce}), we have the
compromise solution for the multi-objective problem
(\ref{eqexstart})-(\ref{eqexend}) as presented in Table 7.
\begin{table}[h]\label{tabresultgc}
\begin{center} {Table 7: Compromise solution using global
criterion method
\begin{tabular}{|l|}
               \hline
               $z_{111} = 2$, $z_{121} = 17$, $z_{211} = 7$, $z_{231} = 24$, $z_{112} = 19$, $z_{122}
               =10$, $z_{232}= 2$,\\ $x_{111}^2 =
64$, $x_{121}^1 = 187$, $x_{121}^2=248$, $x_{211}^1 =79$, $x_{211}^2
=100$, $x_{231}^1 = 312$,\\$x_{231}^2 = 277$, $x_{112}^1 = 261$,
$x_{112}^2=111$, $x_{122}^1 = 173$, $x_{122}^2 = 2$, $x_{232}^1 =
33$,\\$x_{232}^2 = 3$, $Min~f'_{1}=8198.6$, $Min~f'_{2}=771.1$. \\
\hline
\end{tabular} }
\end{center}
\end{table}
\par For multi-objective optimization problem, generally we have to
look for compromise solution (in the sense that there does not exist
a single solution that simultaneously optimizes each of the
objectives) due to the conflicting nature of the objectives. Here
from Table 6 and Table 7, it is observed that if minimization of
transportation cost is given more priority than the transportation
time, then solution in Table 6 is better than the solution in Table
7. Otherwise, if the transportation time is given more priority than
the cost of transportation, then the solution in Table 7 is better.
\begin{figure}\begin{center}
  \includegraphics[width=10cm]{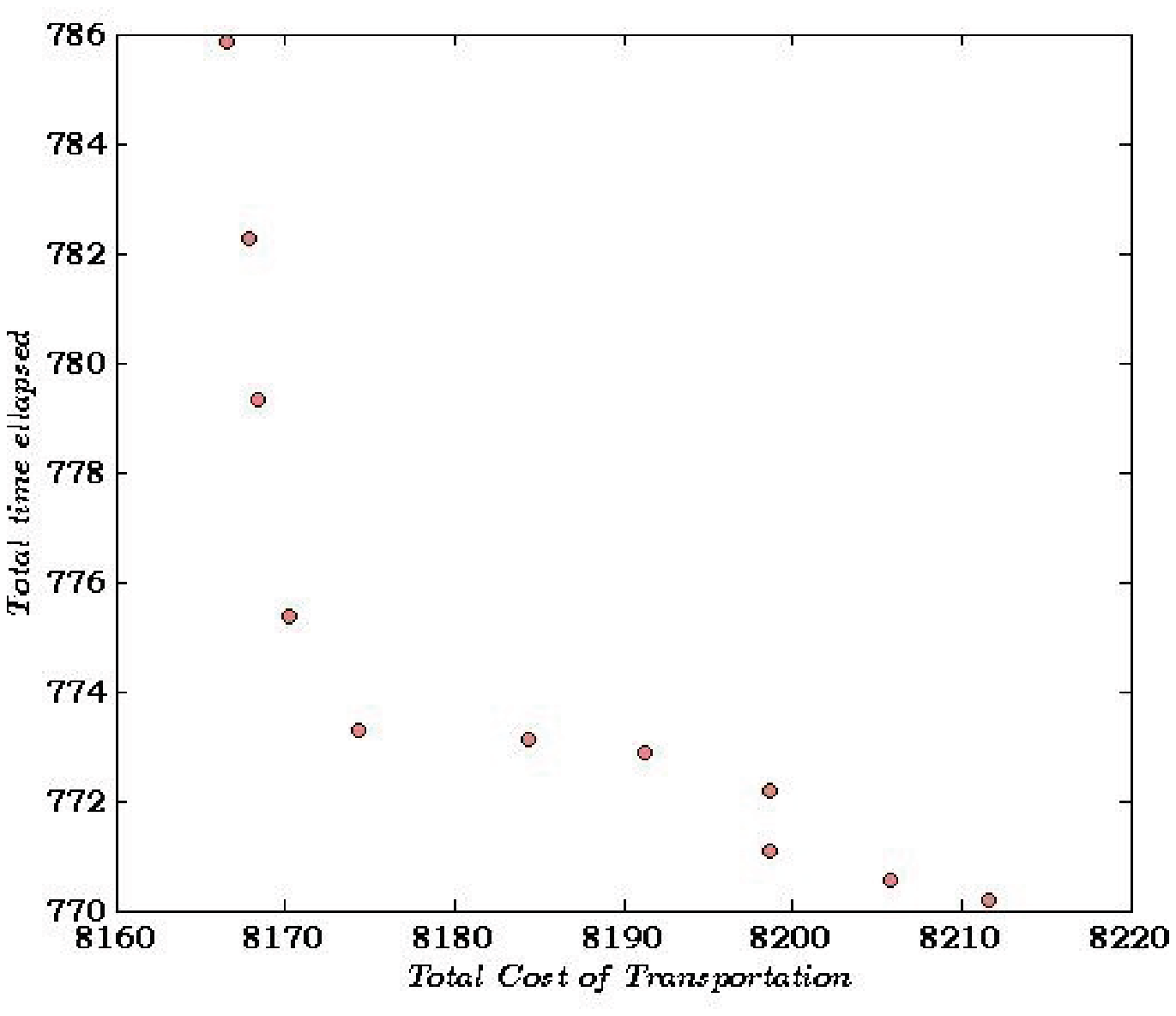}\\
  \caption{Approximated front for the proposed STP.}\label{fig:1}\end{center}
\end{figure}
\par We have also optimize our proposed bi-objective STP by considering several randomly generated weight vectors to obtain multiple non-dominated solutions. Among these solutions we consider the best non-dominated solutions which are depicted in Fig. 1. These solutions are to be consider as members of the approximate front of our proposed bi-objective STP.
\section{Conclusion}
For transportation systems where full vehicles (e.g., light
commercial vehicles, medium duty and heavy duty trucks, rail
coaches, etc.) are to be considered for transportation of products,
different types of issues appear in formulation of the problem. Like
transportation cost of a vehicle which is irrespective of the fact
whether the capacity of the vehicle is filed up or not; volume and
weight capacities of the vehicles; limitation of number of certain
types of vehicles, etc. Also, besides the travel time of a vehicle,
loading and unloading times of products are also important which
depend upon both the vehicle types and product characteristics. In
his paper, we present a multi-objective multi-item solid
transportation model by addressing all those issues. The presented
problem is formulated with transportation time and cost parameters
as fuzzy variables. We have formulated a chance-constrained
programming model to solve the problem with fuzzy parameters, and
then it is transformed into deterministic form. The deterministic
multi-objective problem is solved with two multi-objective
optimization techniques, namely, the fuzzy programming technique and
global criterion method. \par In the presented problem, the
transportation time and cost parameters are considered as usual
(type-1) fuzzy variables, however, this work can be extended with
those parameters as represented by generalized type-2 or interval
type-2 fuzzy variables.

\end{document}